\documentclass[11pt]{amsart}
\usepackage{mathrsfs}

\usepackage{amsmath}
\usepackage{amssymb} %to produce blackboard bold characters(\supsetneq)
\usepackage{enumitem} %lay­out of the three ba­sic list en­vi­ron­ments
\usepackage{mathtools} %an extension package to amsmath + loading amsmath
\usepackage[table]{xcolor} %color of tables
\usepackage[all]{xy} %xymatrix commutative diagram
\usepackage{indentfirst} %Indent first paragraph
\usepackage{setspace} %set­ting the spac­ing be­tween lines

\usepackage[colorlinks,linkcolor=red,anchorcolor=green,citecolor=blue]{hyperref} %color of hyperlinks
\hypersetup{linktocpage = true} %color of hyperlinks of TOC

\usepackage{rotating} %Rotate an image, table or paragraph
\usepackage{longtable} %long table
\newcolumntype{M}[1]{>{\centering\arraybackslash}m{#1}} %define dimension for long stable

\usepackage{mathpazo}
\usepackage{mathrsfs}
\DeclareFontFamily{OMS}{rsfs}{\skewchar\font'60}
\DeclareFontShape{OMS}{rsfs}{m}{n}{<-5>rsfs5 <5-7>rsfs7 <7->rsfs10 }{}
\DeclareSymbolFont{rsfs}{OMS}{rsfs}{m}{n}
\DeclareSymbolFontAlphabet{\scr}{rsfs}
\DeclareSymbolFontAlphabet{\scr}{rsfs}

\usepackage[T1]{fontenc}
\newcommand\sO{{\mathscr O}}

\newcommand{\rk}{{\rm rk}}

\newcommand\bp{{\bar\partial}}

\theoremstyle{plain}
\newtheorem{thm}{Theorem}[section]
\newtheorem{lemma}[thm]{Lemma}
\newtheorem{prop}[thm]{Proposition}
\newtheorem{cor}[thm]{Corollary}
\newtheorem{defn}[thm]{Definition}

\newtheorem{conjecture}[thm]{Conjecture}

\theoremstyle{definition}
\newtheorem{example}[thm]{Example}
\newtheorem{remark}[thm]{Remark}

\newcommand{\btheorem}{\begin{thm}}
	\newcommand{\etheorem}{\end{thm}}
\newcommand{\bproposition}{\begin{prop}}
	\newcommand{\eproposition}{\end{prop}}
\newcommand{\bdefinition}{\begin{defn}}
	\newcommand{\edefinition}{\end{defn}}
\newcommand{\bcorollary}{\begin{cor}}
	\newcommand{\ecorollary}{\end{cor}}
\newcommand{\bproof}{\begin{proof}}
	\newcommand{\eproof}{\end{proof}}
\newcommand{\bremark}{\begin{remark}}
	\newcommand{\eremark}{\end{remark}}
\newcommand{\eexample}{\end{example}}
\newcommand{\bexample}{\begin{example}}

\newcommand{\elemma}{\end{lemma}}
\newcommand{\blemma}{\begin{lemma}}

\newcommand{\sq}{\sqrt{-1}}

\newcommand{\p}{\partial}

\renewcommand{\bar}{\overline}

\renewcommand{\phi}{\varphi}

\newcommand{\beq}{\begin{equation}}
\newcommand{\eeq}{\end{equation}}
\newcommand{\ee}{\end{eqnarray*}}
\newcommand{\be}{\begin{eqnarray*}}

\newcommand{\bd}{\begin{enumerate}}
	\newcommand{\ed}{\end{enumerate}}

\renewcommand{\bp}{\bar{\partial}}

\newcommand{\ts}{\otimes}

\renewcommand{\>}{\rightarrow}

\newcommand{\C}{{\mathbb C}}

\newcommand{\R}{{\mathbb R}}
\numberwithin{equation}{section} %numbering of equations
\makeatletter

\ifnum\@ptsize=0  \fi \ifnum\@ptsize=2  \fi

\usepackage{fancyhdr}
\renewcommand{\leftmark}{{\scriptsize Compact K\"ahler manifolds with quasi-positive second Chern-Ricci curvature}}

\title{Compact K\"ahler manifolds with quasi-positive second Chern-Ricci curvature}

\author{Xiaokui Yang}

\address{Xiaokui Yang, Department of Mathematics and Yau Mathematical Sciences Center, Tsinghua University, Beijing, 100084, China}
\email{xkyang@mail.tsinghua.edu.cn}

\begin{document}
	
	\begin{abstract} Let $X$ be a compact K\"ahler manifold. We prove that if $X$ admits a smooth Hermitian metric $\omega$ with \emph{quasi-positive} second Chern-Ricci curvature $\mathrm{Ric}^{(2)}(\omega)$, then $X$ is projective and rationally connected. In particular, $X$ is simply connected.
	\end{abstract}

	\maketitle
	
	\tableofcontents
	
	\vspace{-0.2cm}
	
	\section{Introduction}
    The geometry of complex manifolds are characterized by various positivity notions in complex differential geometry and algebraic geometry. Since the seminal works of Mori   and Siu-Yau on the solutions to
	Hartshorne
	conjecture   and Frankel conjecture
	(\cite{Mori1979}, \cite{SiuYau1980}) on the characterizations of projective spaces,   many remarkable generalizations  have been
	established,   for instances, Mok's uniformization theorem on compact K\"ahler manifold
	with non-negative holomorphic bisectional curvature (\cite{Mok1988}) and  the  works of Campana,
	Demailly, Peternell and Schneider (\cite{CampanaPeternell1991},
	\cite{DemaillyPeternellSchneider1994})  on the
	structure of projective manifolds with nef tangent bundles. For more geometric characterizations, we refer to \cite{Mok1988,CampanaPeternell1991,DemaillyPeternellSchneider1994,DemaillyPeternellSchneider2001,Wu2002,Chau2012,LiWuZheng2013,Liu2014b,FengLiuWan2017,HuangLeeTamTong2018,NiZheng2019b,Liu2019,LiOuYang2019,LiuOuYang2020} and the references therein.\\
	
	    The  holomorphic sectional curvature also carries much geometric information of   complex manifolds. Indeed, thanks to the breakthrough work \cite{WuYau2016} of Wu-Yau, it is well-known that  a compact K\"ahler manifold $X$ with negative or quasi-negative holomorphic sectional curvature is algebraic and has ample canonical bundle (\cite{TosattiYang2017,WuYau2016a,Diverio2019}), which settles down a long-standing conjecture of S.-T. Yau affirmatively. For more recent works on  non-positive holomorphic sectional curvature, we refer to \cite{Wong1981,Heier2010,Heier2016,Heier2018,Nomura2018,Guenancia2018,YangZheng2019TAMS,WuYau2020} and the references therein. On the other hand, in his “Problem section”,
	    S.-T. Yau proposed the well-known conjecture \cite[Problem~47]{Yau1982} that compact K\"ahler manifolds with positive holomorphic sectional curvature must be projective and rationally connected. Recently, we solved this conjecture affirmatively in \cite{Yang2018} by introducing the concept of RC-positivity for abstract vector bundles, and many properties of such bundles are developed in \cite{Yang2018, Yang2018a,Yang2018b,Yang2021}. For instance, we proved in \cite{Yang2018b}  that a compact K\"ahler manifold with uniformly RC-positive tangent bundle must be projective and rationally connected.  By using the ideas of RC-positivity and some deep analytical techniques in algebraic geoemtry, Shin-ichi Matsumura established in \cite{Matsumura2018a,Matsumura2018b,Matsumura2018} a structure theorem for projective manifolds  with non-negative holomorphic sectional curvature, which is analogous to fundamental works in \cite{Mok1988, Campana2015a,Cao2017a,Cao2019} for manifolds with various non-negative properties (see also an approach in \cite{HeierWong2015}). In the same spirit, Lei Ni and Fangyang Zheng introduced in \cite{Ni2018,NiZheng2018a} various notions of Ricci curvature and scalar curvature to obtain rational connectedness of compact K\"ahler manifolds.\\

	In this paper, we investigate the geometry characterized by the first and second Chern-Ricci curvatures on Hermitian manifolds. Recall that for a Hermitian metric $\omega_g$,  its Chern curvature tensor has components $R_{i\bar j k\bar\ell}$. The first Chern-Ricci curvature is $$\mathrm{Ric}^{(1)}(\omega)=\sq \left(g^{k\bar\ell}R_{i\bar j k\bar \ell}\right)dz^i\wedge d\bar z^j=-\sq \p\bp\log\det(h_{k\bar\ell}) $$
	and the second Chern-Ricci curvature is 
	$$\mathrm{Ric}^{(2)}(\omega)=\sq \left(g^{k\bar\ell}R_{ k\bar \ell i\bar j}\right)dz^i\wedge d\bar z^j.$$
	When the Hermitian metric $\omega_g$ is not K\"ahler,  $\mathrm{Ric}^{(1)}(\omega)$ and $\mathrm{Ric}^{(2)}(\omega)$ are not necessarily the same. It is well-known that $\mathrm{Ric}^{(1)}(\omega)$ represents the first Chern class of the complex manifold. However, the geometry of $\mathrm{Ric}^{(2)}(\omega)$ is still mysterious.\\
	
	Thanks to the celebrated Calabi-Yau theorem (\cite{Yau1978}),  we know that  a
	compact K\"ahler manifold $X$ has a Hermitian metric with positive first Chern-Ricci curvature $\mathrm{Ric}^{(1)}(\omega)$
	if and only if $X$ is Fano. As an analog, we proved in \cite{Yang2018} that if a compact K\"ahler manifold
	admits a smooth Hermitian metric with positive second Chern-Ricci curvature $\mathrm{Ric}^{(2)}(\omega)$, then $X$ is projective and rationally connected. This result is also a generalization of the classical result of Campana \cite{Campana1992} and Koll\'ar-Miyaoka-Mori \cite{KollarMiyaokaMori1992} that Fano manifolds are rationally connected. The main result of this paper is the following theorem.

		\btheorem\label{thm:quasipositivesecondchernriccirationallyconnected} Let $X$ be a compact K\"ahler manifold. If there exist a smooth Hermitian metric $\omega$ on $X$ and a smooth Hermitian metric $h$ on the holomorphic tangent bundle $TX$ such that  $$\mathrm{tr}_\omega R^{(TX,h)}\in \Gamma(X, \mathrm{End}(TX))$$ is quasi-positive. Then $X$ is projective and rationally connected. In particular, $X$ is simply connected.
	\etheorem
\noindent Here quasi-positive means non-negative everywhere and strictly positive at some point.  We follow the ideas in \cite{Campana2015a} and \cite{Graber2003} in the proof of Theorem \ref{thm:quasipositivesecondchernriccirationallyconnected}, and the key new ingredient is an integration argument  for singular Hermitian metrics instead of the pointwise maximum principle for RC-positive vector bundles employed in \cite{Yang2018}, since the latter does not work for manifolds with quasi-positive curvature tensors. It is easy to see that many compact K\"ahler manifolds with stable tangent bundle and positive slope can support smooth Hermitian metrics with positive or quasi-positive second Chern-Ricci curvature as required in Theorem \ref{thm:quasipositivesecondchernriccirationallyconnected} (e.g. \cite{UhlenbeckYau1986,Donaldson1987}). On the other hand, the K\"ahler condition in Theorem \ref{thm:quasipositivesecondchernriccirationallyconnected} is necessary (\cite[Section~6]{LiuYang2017}).
As a special case of Theorem \ref{thm:quasipositivesecondchernriccirationallyconnected}, one has

	\bcorollary \label{quasipositivesecondchernriccirationallyconnected} Let $X$ be a compact K\"ahler manifold. If $X$
	admits a smooth Hermitian metric $\omega$ with quasi-positive second Chern-Ricci curvature $\mathrm{Ric}^{(2)}(\omega)$, then $X$ is projective and rationally connected. In particular, $X$ is simply connected.
	\ecorollary

\noindent	By using  Corollary \ref{quasipositivesecondchernriccirationallyconnected}  and Yau's theorem \cite{Yau1978}, one has the following generalization of the result of  Campana \cite{Campana1992} and Koll\'ar-Miyaoka-Mori \cite{KollarMiyaokaMori1992}.
	
		\bcorollary \label{weakfano} Let $X$ be a compact K\"ahler manifold. If $X$
	admits a smooth Hermitian metric $\omega$ with quasi-positive first Chern-Ricci curvature $\mathrm{Ric}^{(1)}(\omega)$, then $X$ is projective and rationally connected. In particular, $X$ is simply connected.
	\ecorollary

	  \noindent Similarly,	for quasi-negative first Chern-Ricci curvature $\mathrm{Ric^{(1)}}(\omega)$, one has

	\btheorem \label{quasinegative} Let $(X,\omega)$ be a compact Hermitian manifold with quasi-negative first Chern-Ricci curvature $\mathrm{Ric^{(1)}}(\omega)$.
	If $X$ contains no rational curve, then $X$ is  projective  and $K_X$ is ample.
	\etheorem
	
		\noindent  This result is a straightforward consequence of deep results in complex analytical and algebraic geometry (\cite{Mori1982,Siu1984,Demailly1985,Kawamata1985,Birkar2010,Cascini2013}), which would be of independent interest from the viewpoint of complex differential geometry. It is known that  compact complex manifolds with quasi-negative (or quasi-positive) first Chern-Ricci curvature
	$\mathrm{Ric^{(1)}}(\omega)$ are Moishezon (\cite{Siu1984,Demailly1985}), which are not necessarily projective (e.g. \cite{Ma2007}).
	There are also many compact complex manifolds containing no rational curves, for instances,
	hyperbolic manifolds and Hermitian manifolds with non-positive holomorphic sectional curvature. Theorem \ref{quasinegative} also provides another proof of \cite[Corollary~1.1]{Lee2018} that  compact Hermitian manifolds with non-positive holomorphic bisectional curvature and quasi-negative first Chern-Ricci curvature are projective manifolds with ample canonical bundles,  which was established by purely analytical method.\\

	\noindent 
	The following conjecture on quasi-positive holomorphic sectional curvature is well-known and still widely open (e.g. \cite[Conjecture~1.9]{Yang2018b}) and in the special case when $X$ is projective it was confirmed affirmatively by Matsumura  (\cite{Matsumura2018b}). 
	
	\begin{conjecture}\label{quasipositiveHSC} Let $(X,\omega)$ be a compact K\"ahler manifold. If it has quasi-positive holomorphic sectional curvature, then $X$ is projective and rationally connected.
	\end{conjecture}

	\noindent{\bf Acknowledgements.}  The author would like to thank Jie Liu, Wenhao Ou, Valentino Tosatti and S.-T. Yau for inspiring discussions and useful communications.

	\vskip 2\baselineskip

	\section{Singular Hermitian metrics on vector bundles}
	
	Let $X$ be a complex manifold and $\omega$ be a smooth Hermitian metric on $X$. Locally, we can write the curvature tensor of $(T_X, \omega)$ as 
	$$R_{i\bar j k\bar \ell}=-\frac{\p^2 g_{k\bar\ell}}{\p z^i\p\bar z^j}+g^{p\bar q}\frac{\p g_{k\bar q}}{\p z^i}\frac{\p g_{p\bar \ell}}{\p\bar z^j}.$$ where $\omega=\sq g_{i\bar j}dz^i\wedge d\bar z^j$.
	The \emph{first Chern-Ricci curvature} $\mathrm{Ric}^{(1)}(\omega)=\sq R_{i\bar j}dz^i\wedge d\bar z^j$  has components
	$R^{(1)}_{i\bar j}=g^{k\bar\ell}R_{i\bar j k\bar \ell}=-\frac{\p^2\log \det (g_{k\bar\ell})}{\p z^i\p\bar z^j}.$
	The \emph{second Chern-Ricci curvature} is  $\mathrm{Ric}^{(2)}(\omega)=\sq R^{(2)}_{i\bar j}dz^i\wedge d\bar z^j$ where
	$R^{(2)}_{i\bar j}=g^{k\bar\ell}R_{ k\bar \ell i\bar j}$. The scalar curvature $s$ of the Chern connection is $\mathrm{tr}_\omega \mathrm{Ric^{(1)}}(\omega)$ which is also the same as $\mathrm{tr}_\omega \mathrm{Ric^{(2)}}(\omega)$. 
	
		Let $E\>X$ be a holomorphic vector bundle and $h$ be a smooth Hermitian metric on $E$. The curvature $R^E$ of the Chern connection $\nabla $ on $(E,h)$ has a similar formula $$R^E_{i\bar j \alpha\bar\beta}=-\frac{\p^2 h_{\alpha\bar\beta}}{\p z^i\p\bar z^j}+h^{\gamma\bar\delta}\frac{\p h_{\alpha\bar \delta}}{\p z^i}\frac{\p h_{\gamma\bar\beta}}{\p \bar z^j}.$$
	We define 
	$R^{(1)}_{i\bar j}=h^{\alpha\bar\beta} R_{i\bar j \alpha\bar\beta}$ and $R^{(2)}_{\alpha\bar\beta}=g^{i\bar j}R_{i\bar j \alpha\bar\beta}$.
	We also call $\mathrm{tr}_hR^E=\sq R_{i\bar j}^{(1)}dz^i\wedge d\bar z^j$ and $\mathrm{tr}_\omega R^E= \sq R_{\alpha\bar \beta}^{(2)}e^\alpha \ts \bar e^\beta$ the first and second Chern-Ricci curvature of $(E,h)$ with respect to the Hermitian manifold $(X,\omega)$ respectively. When $(E,h)=(TX,\omega)$, they are exactly the same as those curvatures of $(X,\omega)$.   A smooth Hermitian $(1,1)$-form $A=\sq A_{i\bar j}dz^i\wedge d\bar z^j$ on $X$ is call \emph{quasi-positive}, if $(A_{i\bar j})$ is non-negative everywhere and positive at some point of $X$. Similarly, we can define it for a  tensor $A \in \Gamma(X,\mathrm{End}(E))$.\\
	
	 Singular Hermitian metrics on line bundles are introduced in \cite{Demailly1992a} by Demailly. Let $L$ be a holomorphic vector bundle. A singular metric $h^L$ on $L$ can be written locally as $h^L=e^{-\phi}$ for some $\phi\in L^1_{\mathrm{Loc}}(X,\R)$, and the curvature $R^L=-\sq\p\bp\log h^L$ is defined in the sense of distributions. For singular Hermitian metrics on vector bundles, the definition would be very subtle, and we refer to \cite[Section~2]{Paun2018a} for a detailed discussion (see also \cite{Cataldo1998,Raufi2015,DengWangZhangZhou2020}). Recall that, for a smooth Hermitian metric $h$ on $E$, its curvature also takes the local form $R^E=\bp(h^{-1}\p h)$. For the specified purpose in this paper, we only consider singular metrics such that $h^{-1}\p h$ is locally integrable, and then we can use standard theory on distributions (e.g. \cite{Demailly2012}) to define the notion of weak positivity. 
	 
	 \bdefinition \label{definitionforsingularmetric} Let $X$ be a  complex manifold. A vector bundle $E$ is called to have positive second Chern-Ricci curvature in the sense of distributions, if there exist a smooth metric $\omega$ on $X$ and a singular Hermitian metric $h^E$ on $E$ such that
	 \beq \mathrm{tr}_\omega R^E \in \Gamma(X,\mathrm{End}(E))\eeq 
	 is strictly positive in the sense of distributions. The definition for non-negativity, quasi-positivity and etc. can be defined similarly.
	 \edefinition
	 
	 \bremark Of course, the notion of second Chern-Ricci curvature can be defined in a broader way by using similar constructions as in \cite{Cataldo1998}.  So far, it is not clear to the author whether the notions of Griffiths positivity or Nakano positivity for singular metrics defined in  \cite{Paun2018a,Raufi2015,DengWangZhangZhou2020} can imply the positivity of the second Chern-Ricci curvature in suitable sense, though it is obvious for smooth Hermitian metrics.
	\eremark
	
	\vskip 2\baselineskip
	
	\section{Vanishing theorems for singular Hermitian metrics}
	
	The main result of this section is the following theorem.

	\btheorem\label{vanishingforquasipositive} Let $E$ be a holomorphic vector bundle over a compact complex manifold $X$.
	Suppose there exist a smooth Hermitian metric $\omega$ on $X$ and a smooth Hermitian metric $h$ on  $E$ such that  $$\mathrm{tr}_\omega R^{(E,h)}\in \Gamma(X, \mathrm{End}(TX))$$
	is quasi-positive. We have the following assertions.
	
	\bd \item Any invertible \textbf{subsheaf} $L$ of $\sO(\ts^k E^*)$ $(k\geq 1)$ is not pseudo-effective.
	\item $\det E^*$ is not pseudo-effective. \ed
	\etheorem

	\noindent The key difficulty in the proof of Theorem \ref{vanishingforquasipositive} is to deal with a line bundle $L$ which is only a subsheaf of $\sO(\ts^k E^*)$. If it is indeed a subbundle, then its follows from a simple observation (c.f. \cite{Campana2015a}).

	\blemma\label{monotonicityofsecondchernricci} The second Chern-Ricci curvature is decreasing in subbundles and increasing in quotient bundles.
	\elemma

	\noindent It is well-known that the first Chern-Ricci curvature is not necessarily monotone as described in Lemma \ref{monotonicityofsecondchernricci}.
	The proof follows from a standard computation  and we include details here since we need it for the distribution case.
	
	\bproof  Let $(E,h)$ be a Hermitian vector bundle and $S$ be a holomorphic subbundle of $E$.
	Let $r$ be the rank of $E$ and $s$ the rank of $S$. Without
	loss of generality, we can assume,  at a fixed point $p\in X$, there
	exists a local holomorphic frame $\{e_1,\cdots,e_r\}$ of $E$
	centered at point $p$ such that $\{e_1,\cdots, e_s\}$ is a local
	holomorphic frame of $S$. Moreover, we can assume that
	$h(e_\alpha,e_\beta)(p)=\delta_{\alpha\beta}, \text{for } 1\leq
	\alpha, \beta \leq r.$ Hence, the curvature tensor of
	$S$ at point
	$p$ is $$ R^S_{i\bar j \alpha\bar\beta}=-\frac{\p^2
		h_{\alpha\bar\beta}}{\p z^i\p\bar z^j}+\sum_{\gamma=1}^{s}\frac{\p
		h_{\alpha\bar\gamma}}{\p z^i}\frac{\p h_{\gamma\bar\beta}}{\p\bar
		z^j}$$ where $1\leq \alpha,\beta\leq s$. For any Hermitian metric $\omega=\sq g_{i\bar j}dz^i\wedge d\bar z^j$ on $X$, we have
	\beq (\mathrm{tr}_\omega R^E)|_{S}-\mathrm{tr}_\omega R^S= \mathrm{tr}_\omega(R^E|_{S})-\mathrm{tr}_\omega R^S=\sq
	\sum_{\alpha,\beta=1}^sg^{i\bar j}\left(\sum_{\gamma=s+1}^r\frac{\p
		h_{\alpha\bar\gamma}}{\p z^i}\frac{\p h_{\gamma\bar\beta}}{\p\bar
		z^j}\right) e^\alpha\ts
	\bar e^\beta. \label{secondchernriccimonotone}\eeq
	It is easy to see that the right hand side of (\ref{secondchernriccimonotone}) is non-negative. The proof for quotient bundles is similar.
	\eproof
	
	\bcorollary\label{functorialforsecondricci} Let $E$ be a holomorphic vector bundle over a complex manifold $X$.
	
	\bd \item  If $E$ has  positive (resp, non-negative) second Chern-Ricci curvature, then each quotient bundle has    positive (resp, non-negative) second Chern-Ricci curvature.
	
	\item If $E$ has negative (resp, non-positive) second Chern-Ricci curvature, then each subbundle has   negative (resp, non-positive) second Chern-Ricci curvature.
	
	\item If $E$ has  positive (resp, non-negative, quasi-positive,  negative, non-positive, quasi-negative) second Chern-Ricci curvature, then so is $\ts^k E$ for each $k\geq 1$.
	
	\item  If $E$ has  positive (resp, non-negative, quasi-positive,  negative, non-positive, quasi-negative) second Chern-Ricci curvature, then so is $\mathrm{Sym}^{\ts k} E$ $(k\geq 1)$ and $\Lambda^pE$ $ (1\leq p\leq \rk(E))$.

	\ed
	\ecorollary
	
	\bproof We only need to prove $(3)$. From the expression of the induced curvature formula of
	$(\ts ^kE^{},\ts ^kh^{})$, one has
	$$R^{(\ts ^kE^{},\ts ^kh^{})}=\ts ^k R^{(E,h)}\in \Gamma\left(X,\Lambda^{1,1}T^*X\ts \mathrm{End}(\ts^kE)\right),$$
	and
	$$\mathrm{tr}_\omega R^{(\ts ^kE^{},\ts ^kh^{})}=\mathrm{tr}_\omega\left(\ts ^k R^{(E,h)}\right)=\ts^k\left(\mathrm{tr}_\omega R^{(E,h)}\right)\in \Gamma\left(X,\mathrm{End}(\ts^k E)\right).$$
	Hence, the result follows.
	\eproof

	\bremark By standard distribution theory, one has similar results as in   Lemma \ref{monotonicityofsecondchernricci}
and Corollary \ref{functorialforsecondricci}  for singular Hermitian metrics  in  Definition \ref{definitionforsingularmetric}. 	
	\eremark

	\blemma\label{nonwherevanishing} Suppose that $E$ has non-negative second Chern-Ricci curvature in the  sense of distributions and $L$ is a pseudo-effective line bundle.
	If $\sigma\in H^0(X,E^*\ts L^*)$ is a non-trivial holomorphic section, then $\sigma$ does not vanish everywhere.
	\elemma
	
	\bproof Let $\omega$ be a smooth Hermitian metric on $X$ and $h^E$ be a singular metric such that $\mathrm{tr}_\omega R^E$ is non-negative in the sense of distributions.
	Since $$\mathrm{tr}_{e^f\omega} R^{E}=e^{-f}\mathrm{tr}_\omega R^E,$$
	we can assume further that $\omega$ is a Gauduchon metric (\cite{Gauduchon1977}), i.e. $\p\bp\omega^{n-1}=0$
	where $\dim X=n$. Let $h^L$ be a singular metric on $L$ such that its curvature $\theta=-\sq\p\bp\log h^L\geq 0$
	in the sense of distributions. For any $\sigma\in H^0(X,E^*\ts L^*)$, we have
	\beq\mathrm{tr}_\omega\left(\sq \p\bp |\sigma|^2_{h^{E^*}\ts h^{L^*}}\right)
	=|\nabla \sigma|^2_{\omega,h^{E^*},h^{L^*}}+(\mathrm{tr}_\omega R^E\cdot h^{L^*}+\mathrm{tr}_\omega\theta\cdot h^{E^*})(\sigma,\sigma).\label{bochnerformula}\eeq
	Note that we are dealing with positive currents, and the product makes sense.
	An integration by part argument with respect to the Gauduchon metric $\omega$  shows that $\nabla\sigma=0$ a.e., and
	\beq (\mathrm{tr}_\omega\sq\p\bp) |\sigma|^2_{h^{E^*}\ts h^{L^*}}\geq 0 \eeq
	in the sense of distributions.  We know $|\sigma|^2_{h^{E^*}\ts h^{L^*}}$ is a global constant on the compact base $X$.
	If $\sigma$ is non-trivial, this constant is nonzero, i.e. $\sigma$ does not vanish everywhere.
	\eproof
	
	\bcorollary\label{subsheafsubbundle} Suppose that $E$ has non-negative second Chern-Ricci curvature  in the  sense of distributions, and  $L$ is an invertible
	subsheaf of $\sO(E^*)$. If $L$ is pseudo-effective, then $L$ is a line subbundle of $E^*$.
	\ecorollary
	
	\bproof Note that the subsheaf morphism $f: L\> \sO(E^*)$ induces a nonzero section $\sigma \in H^0(X,E^*\ts L^*)$. By Lemma \ref{nonwherevanishing},
	$\sigma$ does not vanish everywhere. This gives a trivial line subbundle  $\underline{\C}$ of the vector bundle $E^*\ts L^*$, and so $L$ is a line subbundle of $E^*$.
	\eproof

	\btheorem\label{distributionvanishing} Suppose that $E$ has quasi-positive second Chern-Ricci curvature  in the  sense of distributions, and  $L$ is a pseudo-effective line bundle. Then
	\beq H^0(X,E^*\ts L^*)=0.\eeq
	\etheorem
	
	\bproof By using a similar argument as in the proof of Lemma \ref{nonwherevanishing}, for any holomorphic section $\sigma\in H^0(X,E^*\ts L^*)$,
	we deduce  from (\ref{bochnerformula}) that the holomorphic section $\sigma$ vanishes on an open subset of $X$. By Aronszajn's principle (\cite{Aronszajn1957}), $\sigma$ is a zero section.
	\eproof

\begin{proof}[Proof of Theorem \ref{vanishingforquasipositive}]
	Let $L$ be an invertible subsheaf of $\sO(\ts^k E^*)$ for some
	$k\geq 1$. Suppose-to the contrary-that $L$ is pseudo-effective. By
	Corollary \ref{functorialforsecondricci}, $\ts^k E$ has
	quasi-positive second Chern-Ricci curvature. By Corollary
	\ref{subsheafsubbundle}, the pseudo-effective invertible subsheaf
	$L$ is actually a line subbundle of $\ts^k E^*$. It is easy to see
	that if $\ts^k E$ has quasi-positive second Chern-Ricci curvature,
	then its dual bundle  $\ts^k E^*$ has quasi-negative second
	Chern-Ricci curvature. By decreasing principle in Lemma
	\ref{monotonicityofsecondchernricci}, the second Chern-Ricci
	curvature of $L$ is quasi-negative. Since $L$ is a line bundle,
	that means, the induced metric $h^L$  has quasi-negative Chern
	scalar curvature $s=\mathrm{tr}_\omega (-\sq\p\bp\log h^L)$. Without
	loss of generality, we can assume $\omega$ is Gauduchon. Therefore \beq
	\int_X s\omega^n=n \int_X c_1^{\mathrm{BC}}(L)\wedge
	\omega^{n-1}<0.\eeq  Hence, $L$ can not be
	pseudo-effective (\cite[Proposition~3.1]{Yang2019CM}, ). This is a contradiction.
	
	On the other hand, since $\det E^*$ is a subbundle of $\ts^k E^*$ for some large $k$, we deduce $\det E^*$ can not be pseudo-effective. \end{proof}

	\vskip 2\baselineskip

	\section{The proof of main theorems}

In this section, we prove Theorem \ref{thm:quasipositivesecondchernriccirationallyconnected}, Corollary \ref{weakfano} and Theorem \ref{quasinegative}.
	
	\bproof[{Proof of Theorem \ref{thm:quasipositivesecondchernriccirationallyconnected}}]  We set $(E,h)=(TX,h)$.  By Corollary \ref{functorialforsecondricci} and Theorem \ref{distributionvanishing}, we have
	$$H^0(X,(T^*X)^{\ts k})=0.$$
	It is well-known that for any $p$ satisfying $1\leq p\leq \dim X$, $\Lambda^pT^*X$ is a direct summand of $(T^*X)^{\ts k}$ for some large $k$. Therefore, we have
	$$H_{\bp}^{p,0}(X,\C)\cong H^{0}(X,\Lambda^p T^*X)=0.$$
	In particular,  $H_{\bp}^{2,0}(X,\C)=0$. Since $X$ is K\"ahler, we deduce that $X$ is projective (\cite{Kodaira1954}).
	Now we following the ideas in \cite{Graber2003,Peternell2006,Campana2015a} to get the rational connectedness.
	By using Theorem \ref{vanishingforquasipositive}, we conclude that $K_X=\det T^*X$ is not pseudo-effective.
	Thanks to \cite{BoucksomDemaillyPuaunPeternell2013}, $X$ is actually uniruled.  Let $\pi:X\dashrightarrow Z$ be the
	associated MRC fibration of $X$. After possibly resolving the
	singularities of $\pi$ and $Z$, we may assume that $\pi$ is a proper
	morphism and $Z$ is smooth. By \cite[Corollary~1.4]{Graber2003}, it
	follows that the target of the MRC fibration is either a point or a
	positive dimensional variety which is not uniruled. Suppose $X$ is
	not rationally connected, then $\dim Z\geq 1$ and $Z$ is not
	uniruled. By using \cite{BoucksomDemaillyPuaunPeternell2013} again, $K_Z$ is indeed pseudo-effective.
	Since $K_Z=\det(T^*Z)$ is  a direct summand of the vector bundle $(T^*Z)^{\ts k}$ for some large $k$ and $\sO((T^*Z)^{\ts k})$ is a subsheaf of $ \sO\left((T^*X)^{\ts k}\right)$,  we obtain a pseudo-effective invertible subsheaf $K_Z$ of $\sO((T^*X)^{\ts k})$, which contradicts to part $(1)$ of Theorem \ref{vanishingforquasipositive}.
	\eproof

\bproof[Proof of Corollary \ref{weakfano}] By the celebrated Calabi-Yau theorem, there exists a smooth \emph{K\"ahler metric} $\omega_0$ on $X$ such that 
\beq \mathrm{Ric}^{(1)}(\omega_0)=\mathrm{Ric}^{(1)}(\omega).\eeq
Since $\omega_0$ is K\"ahler, one can deduce
$$\mathrm{Ric}^{(2)}(\omega_0)=\mathrm{Ric}^{(1)}(\omega_0).$$
Therefore, Corollary \ref{weakfano} follows from Corollary \ref{quasipositivesecondchernriccirationallyconnected}.
\eproof
	
	\bproof[Proof of Theorem \ref{quasinegative}] By using the criterion of Siu-Demailly (\cite{Siu1984,Demailly1985}), we know $X$ is a Moishezon manifold and $K_X$ is  big. On the other hand, by the deep result  \cite[Corollary~1.4.6]{Birkar2010}, it is proved in \cite[Theorem~3.1]{Cascini2013} that a Moishezon manifold without rational curve must be projective. Moreover, by the base-point-free theorem and the relative cone theorem \cite{Mori1982,Kawamata1985,KollarMiyaokaMori1992}, one can deduce that $K_X$ is ample since $X$ is a projective manifold of general type  without rational curve.
	\eproof

\vskip 2\baselineskip

	\def\cprime{$'$} %newcommand of prime over letters

	\renewcommand\refname{Reference}
	\bibliographystyle{alpha}
	\bibliography{RC2}
	
\end{document}